\documentclass[12pt]{amsart}
\usepackage{amsmath,amsthm}
\numberwithin{equation}{section}
\newtheorem{thm}[equation]{Theorem}
\newtheorem{cor}[equation]{Corollary}

\newtheorem{prp}[equation]{Proposition}

\theoremstyle{definition}
\newtheorem{df}[equation]{Definition}

\newtheorem{conj}[equation]{Conjecture}
\theoremstyle{remark}
\newtheorem{rem}[equation]{Remark}
\newcommand{\sprf}{\noindent{\it Proof.}}
\newcommand{\sqed}{\hfill\rule{1.3mm}{3mm}\medskip}

\newcounter{stareq}
\def\thestareq{\fnsymbol{stareq}}

\newcommand{\bd}{\begin{description}}
\newcommand{\ed}{\end{description}}
\begin{document}

\date{\today}

\title{Further Developments of Sinai's Ideas: \\ The Boltzmann-Sinai Hypothesis}

\author{N\'andor Sim\'anyi}

\address{The University of Alabama at Birmingham \\
Department of Mathematics \\
1300 University Blvd., Suite 452 \\
Birmingham, AL 35294 U.S.A.}

\email{simanyi@uab.edu}



\begin{abstract}
  In 1963 Ya. G. Sinai \cite{Sin(1963)} formulated a modern version of
  Boltzmann's ergodic hypothesis, what we now call the
  ``Boltzmann-Sinai Ergodic Hypothesis'': The billiard system of $N$
  ($N\ge 2$) hard balls of unit mass moving on the flat torus
  $\mathbb{T}^\nu=\mathbb{R}^\nu/\mathbb{Z}^\nu$ ($\nu\ge 2$) is
  ergodic after we make the standard reductions by fixing the values
  of trivial invariant quantities. It took fifty years and the
  efforts of several people, including Sinai himself, until this
  conjecture was finally proved. In this short survey we provide a
  quick review of the closing part of this process, by showing how
  Sinai's original ideas developed further between 2000 and 2013,
  eventually leading the proof of the conjecture.
\end{abstract}

\maketitle

\section{Posing the Problem \\ The Investigated Models} \label{models}

Non-uniformly hyperbolic systems (possibly, with singularities) play a pivotal
role in the ergodic theory of dynamical systems. Their systematic study
started several decades ago, and it is not our goal here to provide the reader
with a comprehensive review of the history of these investigations but,
instead, we opt for presenting in nutshell a cross section of a few selected
results.

In 1939 G. A. Hedlund and E. Hopf \cite{He(1939)}, \cite{Ho(1939)}, proved the
hyperbolic ergodicity of geodesic flows on closed, compact surfaces with
constant negative curvature by inventing the famous method of "Hopf chains"
constituted by local stable and unstable invariant manifolds.

In 1963 Ya. G. Sinai \cite{Sin(1963)} formulated a modern version of Boltzmann's
ergodic hypothesis, what we call now the "Boltzmann-Sinai Hypothesis":
the billiard system of $N$ ($\ge2$) hard balls of unit mass moving on the
flat torus $\mathbb{T}^\nu=\mathbb{R}^\nu/\mathbb{Z}^\nu$ ($\nu\ge2$) is ergodic after
we make the standard reductions by fixing the values of the trivial invariant
quantities. It took seven years until he proved this conjecture for the
case $N=2$, $\nu=2$ in \cite{Sin(1970)}. Another $17$ years later N. I. Chernov
and Ya. G. Sinai \cite{S-Ch(1987)} proved the hypothesis for the case $N=2$,
$\nu\ge 2$ by also proving a powerful and very useful theorem on local
ergodicity.

In the meantime, in 1977, Ya. Pesin \cite{P(1977)} laid down the foundations
of his theory on the ergodic properties of smooth, hyperbolic dynamical
systems. Later on this theory (nowadays called Pesin theory) was
significantly extended by A. Katok and J-M. Strelcyn \cite{K-S(1986)}
to hyperbolic systems with singularities. That theory is already
applicable for billiard systems, too.

Until the end of the seventies the phenomenon of hyperbolicity (exponential
instability of trajectories) was almost exclusively attributed to some
direct geometric scattering effect, like negative curvature of space, or
strict convexity of the scatterers. This explains the profound shock that
was caused by the discovery of L. A. Bunimovich \cite{B(1979)}: Certain focusing
billiard tables (like the celebrated stadium) can also produce complete
hyperbolicity and, in that way, ergodicity. It was partly this result that
led to Wojtkowski's theory of invariant cone fields, \cite{W(1985)},
\cite{W(1986)}.

The big difference between the system of two balls in $\mathbb{T}^\nu$
($\nu\ge 2$, \cite{S-Ch(1987)}) and the system of $N$ ($\ge 3$) balls in
$\mathbb{T}^\nu$ is that the latter one is merely a so called semi-dispersive
billiard system (the scatterers are convex but not strictly convex
sets, namely cylinders), while the former one is strictly dispersive
(the scatterers are strictly convex sets). This fact makes the proof
of ergodicity (mixing properties) much more complicated. In our series
of papers jointly written with A. Kr\'amli and D. Sz\'asz \cite{K-S-Sz(1990)},
\cite{K-S-Sz(1991)}, and \cite{K-S-Sz(1992)}, we managed to prove the (hyperbolic)
ergodicity of three and four billiard balls on the toroidal container
$\mathbb{T}^\nu$. By inventing new topological methods and the Connecting Path
Formula (CPF), in the two-part paper \cite{Sim(1992)-I}, \cite{Sim(1992)-II} I proved the (hyperbolic)
ergodicity of $N$ hard balls in $\mathbb{T}^\nu$, provided that $N\le\nu$.

The common feature of hard ball systems is -- as D. Sz\'asz pointed this
out first in \cite{Sz(1993)} and \cite{Sz(1994)} -- that all of theom belong to the
family of so called cylindric billiards, the definition of which can be
found later in this survey. However, the first appearance of a special,
3-D cylindric billiard system took place in \cite{K-S-Sz(1989)}, where we
proved the ergodicity of a 3-D billiard flow with two orthogonal
cylindric scatterers. Later D. Sz\'asz \cite{Sz(1994)} presented a complete
picture (as far as ergodicity is concerned) of cylindric billiards with
cylinders whose generator subspaces are spanned by mutually orthogonal
coordinate axes. The task of proving ergodicity for the first non-trivial,
non-orthogonal cylindric billiard system was taken up in \cite{S-Sz(1994)}.

Finally, in our joint venture with D. Sz\'asz \cite{S-Sz(1999)} we managed to
prove the complete hyperbolicity of {\emph typical} hard ball systems on flat tori.

\subsection{Cylindric billiards} Consider the $d$-dimensional
($d\ge 2$) flat torus $\mathbb{T}^d=\mathbb{R}^d/\mathcal{L}$ supplied with the
usual Riemannian inner product $\langle\, .\, ,\, .\, \rangle$ inherited
from the standard inner product of the universal covering space $\mathbb{R}^d$.
Here $\mathcal{L}\subset\mathbb{R}^d$ is supposed to be a lattice, i. e. a discrete
subgroup of the additive group $\mathbb{R}^d$ with $\rm{rank}(\mathcal{L})=d$.
The reason why we want to allow general lattices other than just the
integer lattice $\mathbb{Z}^d$ is that otherwise the hard ball systems would
not be covered. The geometry of the structure lattice $\mathcal{L}$ in the
case of a hard ball system is significantly different from the geometry
of the standard orthogonal lattice $\mathbb{Z}^d$ in the Euclidean space
$\mathbb{R}^d$.

The configuration space of a cylindric billiard is
$\mathbf{Q}=\mathbb{T}^d\setminus\left(C_1\cup\dots\cup C_k\right)$, where the
cylindric scatterers $C_i$ ($i=1,\dots,k$) is defined as follows:

Let $A_i\subset\mathbb{R}^d$ be a so called lattice subspace of $\mathbb{R}^d$,
which means that $\rm{rank}(A_i\cap\mathcal{L})=\rm{dim}A_i$. In this case
the factor $A_i/(A_i\cap\mathcal{L})$ is a subtorus in $\mathbb{T}^d=\mathbb{R}^d/\mathcal{L}$,
which will be taken as the generator of the cylinder 
$C_i\subset\mathbb{T}^d$, $i=1,\dots,k$. Denote by $L_i=A_i^\perp$ the
orthocomplement of $A_i$ in $\mathbb{R}^d$. Throughout this survey article we will
always assume that $\rm{dim}L_i\ge 2$. Let, furthermore, the numbers
$r_i>0$ (the radii of the spherical cylinders $C_i$) and some translation
vectors $t_i\in\mathbb{T}^d=\mathbb{R}^d/\mathcal{L}$ be given. The translation
vectors $t_i$ play a crucial role in positioning the cylinders $C_i$
in the ambient torus $\mathbb{T}^d$. Set
\[
C_i=\left\{x\in\mathbb{T}^d:\; \text{dist}\left(x-t_i,\, A_i/(A_i\cap\mathcal{L})\right)<r_i \right\}.
\]
In order to avoid further unnecessary complications, we always assume that
the interior of the configuration space 
$\mathbf{Q}=\mathbb{T}^d\setminus\left(C_1\cup\dots\cup C_k\right)$ is connected.
The phase space $\mathbf{M}$ of our cylindric billiard flow will be the
unit tangent bundle of $\mathbf{Q}$ (modulo some natural gluings at its
boundary), i. e. $\mathbf{M}=\mathbf{Q}\times\mathbb{S}^{d-1}$. (Here $\mathbb{S}^{d-1}$
denotes the unit sphere of $\mathbb{R}^d$.)

The dynamical system $(\mathbf{M},\, \{S^t\},\, \mu)$, where $S^t$ ($t\in\mathbb{R}$) is the dynamics 
defined by uniform motion inside the domain $\mathbf{Q}$ and specular
reflections at its boundary (at the scatterers), and $\mu$ is the
Liouville measure, is called a cylindric billiard flow we investigate.

\subsection{Transitive cylindric billiards}
The main conjecture concerning the (hyperbolic) ergodicity of cylindric
billiards is the "Erd\H otarcsa conjecture" (named after the picturesque
village in rural Hungary where it was initially formulated) that appeared
as Conjecture 1 in Section 3 of \cite{S-Sz(2000)}:

\begin{conj}{The Erd\H otarcsa Conjecture} A cylindric billiard flow is
  ergodic if and only if it is transitive, i.e. the Lie group
  generated by all rotations across the constituent spaces of the
  cylinders acts transitively on the sphere of compound velocities,
  see Section 3 of \cite{S-Sz(2000)}. In the case of transitivity the
  cylindric billiard system is actually a completely hyperbolic
  Bernoulli flow, see \cite{C-H(1996)} and \cite{O-W(1998)}.
\end{conj}

The theorem of \cite{Sim(2002)} proves a slightly relaxed version of this
conjecture (only full hyperbolicity without ergodicity) for a wide class
of cylindric billiard systems, namely the so called "transverse systems",
which include every hard ball system.

\subsection{Transitivity} Let $L_1,\dots,L_k\subset\mathbb{R}^d$
be subspaces, $A_i=L_i^\perp$, $\text{dim}L_i\ge 2$, $i=1,\dots,k$. Set
\[
\mathcal{G}_i=\left\{U\in\rm{SO}(d):\, U\big|A_i=\rm{Id}_{A_i}\right\},
\]
and let $\mathcal{G}=\left\langle\mathcal{G}_1,\dots,\mathcal{G}_k\right\rangle\subset\rm{SO}(d)$
be the algebraic generate of the compact, connected Lie subgroups
$\mathcal{G}_i$ in $\rm{SO}(d)$. The following notions appeared in Section 3 of
\cite{S-Sz(2000)}.

\begin{df} We say that the system of base spaces
$\{L_1,\dots,L_k\}$ (or, equivalently, the cylindric billiard system defined
by them) is {\emph transitive} if and only if the group $\mathcal{G}$ acts 
transitively on the unit sphere $\mathbb{S}^{d-1}$ of $\mathbb{R}^d$.
\end{df}

\begin{df} We say that the system of subspaces
$\{L_1,\dots,L_k\}$ has the Orthogonal Non-splitting Property (ONSP) if there
is no non-trivial orthogonal splitting $\mathbb{R}^d=B_1\oplus B_2$ of
$\mathbb{R}^d$ with the property that for every index $i$ ($1\le i\le k$)
$L_i\subset B_1$ or $L_i\subset B_2$.
\end{df}

The next result can be found in Section 3 of \cite{S-Sz(2000)} (see 3.1--3.6
thereof):

\begin{prp} For the system of subspaces
$\{L_1,\dots,L_k\}$ the following three properties are equivalent:

(1) $\{L_1,\dots,L_k\}$ is transitive;

(2) $\{L_1,\dots,L_k\}$ has the ONSP;

(3) the natural representation of $\mathcal{G}$ in $\mathbb{R}^d$ is irreducible.
\end{prp}
  
\subsection{Transverseness}

\begin{df} We say that the system of subspaces
$\{L_1,\dots,L_k\}$ of $\mathbb{R}^d$ is {\emph transverse} if the following
property holds: For every {\emph non-transitive} subsystem
$\{L_i:\, i\in I\}$ ($I\subset\{1,\dots,k\}$) there exists an index 
$j_0\in\{1,\dots,k\}$ such that $P_{E^+}(A_{j_0})=E^+$, where
$A_{j_0}=L^\perp_{j_0}$, and $E^+=\text{span}\{L_i:\, i\in I\}$. We note
that in this case, necessarily, $j_0\not\in I$, otherwise $P_{E^+}(A_{j_0})$
would be orthogonal to the subspace $L_{j_0}\subset E^+$. Therefore, every
transverse system is automatically transitive.
\end{df}

We note that every hard ball system is transverse, see \cite{Sim(2002)}. The main result of the paper is the following theorem.

\begin{thm} Assume that the cylindric billiard system is
transverse. Then this billiard flow is completely
hyperbolic, i.e. all relevant Lyapunov exponents are nonzero almost
everywhere. Consequently, such dynamical systems have (at most countably
many) ergodic components of positive measure, and the restriction of
the flow to the ergodic components has the Bernoulli property,
see \cite{C-H(1996)} and \cite{O-W(1998)}.
\end{thm}

Am immediate consequence of this result is

\begin{cor}
  Every hard ball system is completely hyperbolic.
\end{cor}

Thus, the theorem of \cite{Sim(2002)} generalizes the main result of 
\cite{S-Sz(1999)}, where the complete hyperbolicity of {\emph almost every}
hard ball system was proven.

\bigskip

\section{Toward Ergodicity} \label{ergodicity}

In the series of articles \cite{K-S-Sz(1989)}, \cite{K-S-Sz(1991)}, 
\cite{K-S-Sz(1992)}, \cite{Sim(1992)-I}, and \cite{Sim(1992)-II}
the authors developed a powerful, three-step strategy for 
proving the (hyperbolic) ergodicity of hard ball systems. First of all,
all these proofs are inductions on the number $N$ of balls involved in the
problem. Secondly, the induction step itself consists of the following three
major steps:

\subsection{Step I} To prove that every non-singular (i. e. smooth)
trajectory segment $S^{[a,b]}x_0$ with a ``combinatorially rich'' symbolic collision sequence is automatically sufficient
(or, in other words, ``geometrically hyperbolic''), provided that the phase point $x_0$ does not belong to a countable union $J$
of smooth sub-manifolds with codimension at least two. (Containing the 
exceptional phase points.)

Here combinatorial richness means that the symbolic collision sequence of the orbit segment contains a large enough number
of consecutive, connected collision graphs, see also the introductory section of \cite{S-Sz(1999)}.

The exceptional set $J$ featuring this result is negligible in our dynamical
considerations -- it is a so called slim set, i.e. a subset of the phase space
$\mathbf{M}$ that can be covered by a countable union $\bigcup_{n=1}^\infty F_n$ of closed, zero-measured
subsets $F_n$ of $\mathbf{M}$ that have topological co-dimension at least $2$.

\subsection{Step II} Assume the induction hypothesis, i. e. that all hard
ball systems with $N'$ balls ($2\le N'<N$) are (hyperbolic and) ergodic. 
Prove that then there exists a slim set $S\subset\mathbf{M}$ with the following property:
For every phase point $x_0\in\mathbf{M}\setminus S$ the whole trajectory
$S^{(-\infty,\infty)}x_0$ contains at most one singularity and its symbolic collision
sequence is combinatorially rich, just as required by the result of Step I.

\subsection{Step III} By using again the induction hypothesis, prove that
almost every singular trajectory is sufficient in the time interval
$(t_0,\, \infty)$, where $t_0$ is the time moment of the singular reflection.
(Here the phrase ``almost every'' refers to the volume defined by the induced
Riemannian metric on the singularity manifolds.)

We note here that the almost sure sufficiency of the singular trajectories
(featuring Step III) is an essential condition for the proof of the celebrated
Theorem on Local Ergodicity for algebraic semi-dispersive billiards 
proved by Chernov and Sinai in \cite{S-Ch(1987)}.
Under this assumption the theorem of \cite{S-Ch(1987)} states that in any 
algebraic semi-dispersive billiard system (i.e. in a system such that the smooth
components of the boundary $\partial\mathbf{Q}$ are algebraic hypersurfaces)
a suitable, open neighborhood $U_0$ of any hyperbolic phase point 
$x_0\in\mathbf{M}$ (with at most one singularity on its trajectory) belongs
to a single ergodic component of the billiard flow. 

In an inductive proof of ergodicity, steps I and II together ensure that
there exists an arc-wise connected set
$C\subset\mathbf{M}$ with full measure, such that every phase point $x_0\in C$
is hyperbolic with at most one singularity on its trajectory. Then the cited
Theorem on Local Ergodicity (now taking advantage of the result of Step III)
states that for every phase point $x_0\in C$ an open neighborhood $U_0$ of
$x_0$ belongs to one ergodic component of the flow. Finally, the connectedness
of the set $C$ and $\mu(C)=1$ easily imply that the billiard flow with $N$ balls
is indeed ergodic, and actually fully hyperbolic, as well.

In the papers \cite{S-Sz(1999)}, \cite{Sim(2003)}, and \cite{Sim(2004)} we investigated systems
of hard balls with masses $m_1,m_2,\dots,m_N$ ($m_i>0$) moving on the flat torus
$\mathbb{T}_L^\nu=\mathbb{R}^\nu/L\cdot \mathbb{Z}^\nu$, $L>0$.

The main results of the papers \cite{Sim(2003)} and \cite{Sim(2004)} are summarized as follows:

\begin{thm} For almost every selection $(m_1,\dots,m_N;\, L)$ of the external geometric
parameters from the region $m_i>0$, $L>L_0(r,\,\nu)$, where the interior 
of the phase space is connected, it is true that the billiard flow
$\left(\mathbf{M}_{\vec m,L},\{S^t\},\mu_{\vec m,L}\right)$ of the $N$-ball
system is ergodic and completely hyperbolic. Then, following from the results
of \cite{C-H(1996)} and \cite{O-W(1998)}, such a semi-dispersive billiard system actually enjoys
the Bernoulli mixing property, as well.
\end{thm}

\begin{rem} We note that the results of the papers \cite{Sim(2003)} and \cite{Sim(2004)}
nicely complement each other. They precisely assert the same, almost
sure ergodicity of hard ball systems in the cases $\nu=2$ and
$\nu\ge 3$, respectively. It should be noted, however, that the proof of
\cite{Sim(2003)} is primarily dynamical-geometric (except the verification
of the Chernov-Sinai Ansatz), whereas the novel parts of \cite{Sim(2004)}
are fundamentally algebraic.
\end{rem}

\begin{rem} The above inequality $L>L_0(r,\,\nu)$ corresponds to 
physically relevant situations. Indeed, in the case $L<L_0(r,\,\nu)$ the 
particles would not have enough room to even freely exchange positions. 
\end{rem}

\bigskip

\section{The Conditional Proof} \label{conditional}

In the paper \cite{Sim(2009)} we again considered the system of $N$
($\ge 2$) elastically colliding hard spheres with masses $m_1,\dots,m_N$
and radius $r$ on the flat unit torus $\mathbb{T}^\nu$, $\nu\ge 2$. We
proved the Boltzmann-Sinai Ergodic Hypothesis, i. e. the full
hyperbolicity and ergodicity of such systems for every selection
$(m_1,\dots,m_N;\, r)$ of the external parameters, provided that almost
every singular orbit is geometrically hyperbolic (sufficient),
i. e. the so called Chernov-Sinai Ansatz is true.  The proof
does not use the formerly developed, rather involved algebraic
techniques, instead it extensively employs dynamical methods and tools
from geometric analysis.

To upgrade the full hyperbolicity to ergodicity, one needs to refine
the analysis of the degeneracies, i.e. the set of non-hyperbolic phase
points.  For hyperbolicity, it was enough that the degeneracies made a
subset of codimension $\ge 1$ in the phase space.  For ergodicity, one
has to show that its codimension is $\ge 2$, or to find some other ways
to prove that the (possibly) arising one-codimensional, smooth submanifolds of
non-sufficiency are incapable of separating distinct, open ergodic
components from each other. The latter approach was successfully pursued in \cite{Sim(2009)}. In the
paper \cite{Sim(2003)} I took the first step in the direction of proving
that the codimension of exceptional manifolds is at least two: It was
proved there that the systems of $N\ge 2$ disks on a 2D torus (i.e., $\nu=2$)
are ergodic for typical (generic) $(N+1)$-tuples of external
parameters $(m_1,\dots,m_N,r)$. The proof involved some
algebro-geometric techniques, thus the result is restricted to generic
parameters $(m_1,\dots,m_N;\,r)$.  But there was a good reason to
believe that systems in $\nu\ge 3$ dimensions would be somewhat easier
to handle, at least that was indeed the case in early studies.

In the paper \cite{Sim(2004)} I was able to further improve the
algebro-geometric methods of \cite{S-Sz(1999)}, and proved that for any $N\ge 2$,
$\nu\ge 2$, and for almost every selection $(m_1,\dots,m_N;\,r)$ of the external
geometric parameters the corresponding system of $N$ hard balls on 
$\mathbb{T}^\nu$ is (fully hyperbolic and) ergodic.

In the paper \cite{Sim(2009)} the following result was obtained.

\begin{thm} For any integer values $N\ge 2$, $\nu\ge 2$,
and for every $(N+1)$-tuple $(m_1,\dots,m_N,r)$ of the external geometric
parameters the standard hard ball system 
$\left(\mathbf{M}_{\vec m,r},\,\left\{S_{\vec m,r}^t\right\},\, \mu_{\vec m,r}\right)$
is (fully hyperbolic and) ergodic, provided that the Chernov-Sinai Ansatz holds true
for all such systems.
\end{thm}

\begin{rem} The novelty of the theorem (as compared to the result 
in \cite{Sim(2004)}) is that it applies to every $(N+1)$-tuple of external 
parameters (provided that the interior of the phase space is connected), 
without an exceptional zero-measure set. Somehow, the most annoying shortcoming
of several earlier results was exactly the fact that those results were
only valid for hard sphere systems apart from an undescribed, countable
collection of smooth, proper submanifolds of the parameter space 
$\mathbb{R}^{N+1}\ni(m_1,m_2,\dots,m_N;\,r)$. Furthermore, those proofs do not
provide any effective means to check if a given $(m_1,\dots,m_N;\,r)$-system
is ergodic or not, most notably for the case of equal masses in Sinai's
classical formulation of the problem.
\end{rem}

\begin{rem} The present result speaks about exactly the same 
models as the result of \cite{Sim(2002)}, but the statement of this new theorem is 
obviously stronger than that of the theorem in \cite{Sim(2002)}: It has been known
for a long time that, for the family of semi-dispersive billiards, ergodicity
cannot be obtained without also proving full hyperbolicity.
\end{rem}

\begin{rem} As it follows from the results of \cite{C-H(1996)} and
\cite{O-W(1998)}, all standard hard ball systems (the models covered by the
theorems of this survey), once they are proved to be mixing, they also enjoy the much stronger
Bernoulli mixing property. However, even the K-mixing property 
of semi-dispersive billiard systems follows from their ergodicity, as the
classical results of Sinai in \cite{Sin(1968)}, \cite{Sin(1970)}, and \cite{Sin(1979)} show.
\end{rem}

In the subsequent part of this survey we review the necessary technical
prerequisites of the proof, along with some of the needed references
to the literature. The fundamental objects of the paper \cite{Sim(2009)} are the so
called "exceptional manifolds" or "separating manifolds" $J$: they are
codimension-one submanifolds of the phase space that are separating
distinct, open ergodic components of the billiard flow.

In \S3 of \cite{Sim(2009)} we proved Main Lemma 3.5, which states, roughly speaking,
the following: Every separating manifold $J\subset\mathbf{M}$ contains at
least one sufficient (or geometrically hyperbolic) phase
point. The existence of such a sufficient phase point $x\in J$,
however, contradicts the Theorem on Local Ergodicity of Chernov and
Sinai (Theorem 5 in \cite{S-Ch(1987)}), since an open neighborhood $U$ of $x$
would then belong to a single ergodic component, thus violating the
assumption that $J$ is a separating manifold. In \S4 this result was
exploited to carry out an inductive proof of the (hyperbolic)
ergodicity of every hard ball system, provided that the Chernov-Sinai
Ansatz holds true for all hard ball systems.

In what follows, we make an attempt to briefly outline the key ideas
of the proof of Main Lemma 3.5 of \cite{Sim(2009)}. Of course, this
outline will lack the majority of the nitty-gritty details,
technicalities, that constitute an integral part of the proof. The proof is a proof
by contradiction.

We consider the one-sided, tubular neighborhoods $U_\delta$ of $J$ with
radius $\delta>0$. Throughout the whole proof of the main
lemma the asymptotics of the measures $\mu(X_\delta)$ of certain
(dynamically defined) sets $X_\delta\subset U_\delta$ are studied, as
$\delta\to 0$. We fix a large constant $c_3\gg 1$, and for typical
points $y\in U_\delta\setminus U_{\delta/2}$ (having non-singular
forward orbits and returning to the layer $U_\delta\setminus
U_{\delta/2}$ infinitely many times in the future) we define the
arc-length parametrized curves $\rho_{y,t}(s)$ ($0\le s\le h(y,t)$) in
the following way: $\rho_{y,t}$ emanates from $y$ and it is the curve
inside the manifold $\Sigma_0^t(y)$ with the steepest descent towards
the separating manifold $J$. Here $\Sigma_0^t(y)$ is the inverse image
$S^{-t}\left(\Sigma_t^t(y)\right)$ of the flat, local orthogonal
manifold passing through $y_t=S^t(y)$. The terminal point $\Pi(y)=\rho_{y,t}\left(h(y,t)\right)$
of the smooth curve $\rho_{y,t}$ is either

\medskip

(a) on the separating manifold $J$, or

\medskip

(b) on a singularity of order $k_1=k_1(y)$.

\medskip

\noindent
The case (b) is further split in two sub-cases, as follows:

\medskip

(b/1) $k_1(y)<c_3$;

\medskip

(b/2) $c_3\le k_1(y)<\infty$.

\medskip

About the set $\overline{U}_\delta(\infty)$ of (typical) points $y\in U_\delta\setminus U_{\delta/2}$ with
property (a) it is shown that, actually, $\overline{U}_\delta(\infty)=\emptyset$.
Roughly speaking, the reason for this is the following: For a point
$y\in \overline{U}_\delta(\infty)$ the powers $S^t$ of the flow
exhibit arbitrarily large contractions on the curves $\rho_{y,t}$,
thus the infinitely many returns of $S^t(y)$ to the
layer $U_\delta\setminus U_{\delta/2}$ would "pull up" the other
endpoints $S^t\left(\Pi(y)\right)$ to the region $U_\delta\setminus
J$, consisting entirely of sufficient points, and showing that the
point $\Pi(y)\in J$ itself is sufficient, thus violating the indirect hypothesis.

The set $\overline{U}_\delta\setminus\overline{U}_\delta(c_3)$ of all
phase points $y\in U_\delta\setminus U_{\delta/2}$ with the property
$k_1(y)<c_3$ are dealt with by a lemma, where it is shown that
\[
\mu\left(\overline{U}_\delta\setminus\overline{U}_\delta(c_3)\right)=o(\delta),
\]
as $\delta\to 0$. The reason, in rough terms, is that such phase
points must lie at the distance $\le\delta$ from the compact
singularity set
\[
\bigcup_{0\le t\le 2c_3}S^{-t}\left(\mathcal{S}\mathcal{R}^-\right),
\]
and this compact singularity set is transversal to $J$, thus ensuring the measure estimate
$\mu\left(\overline{U}_\delta\setminus\overline{U}_\delta(c_3)\right)=o(\delta)$.

Finally, the set $F_\delta(c_3)$ of (typical) phase points 
$y\in U_\delta\setminus U_{\delta/2}$ with $c_3\le k_1(y)<\infty$ is
dealt with by lemmas 3.36, 3.37, and Corollary 3.38 of \cite{Sim(2009)}, where it is shown
that $\mu\left(F_\delta(c_3)\right)\le C\cdot \delta$, with constants
$C$ that can be chosen arbitrarily small by selecting the constant
$c_3\gg 1$ big enough. The ultimate reason of this measure estimate is
the following fact: For every point $y\in F_\delta(c_3)$ the
projection
\[
\tilde{\Pi}(y)=S^{t_{\overline{k}_1(y)}}\in\partial\mathbf{M}
\]
(where $t_{\overline{k}_1(y)}$ is the time of the
$\overline{k}_1(y)$-th collision on the forward orbit of $y$) will
have a tubular distance $z_{tub}\left(\tilde{\Pi}(y)\right)\le
C_1\delta$ from the singularity set $\mathcal{S}\mathcal{R}^-\cup\mathcal{S}\mathcal{R}^+$,
where the constant $C_1$ can be made arbitrarily small by choosing the
contraction coefficients of the powers $S^{t_{\overline{k}_1(y)}}$ on
the curves $\rho_{y,t_{\overline{k}_1(y)}}$ arbitrarily small with the
help of the result in Appendix II. The upper masure estimate (inside
the set $\partial\mathbf{M}$) of the set of such points 
$\tilde{\Pi}(y)\in\partial\mathbf{M}$ (Lemma 2 in \cite{S-Ch(1987)}) finally
yields the required upper bound 
$\mu\left(F_\delta(c_3)\right)\le C\cdot \delta$ with arbitrarily
small positive constants $C$ (if $c_3\gg 1$ is big enough).

The listed measure estimates and the obvious fact
\[
\mu\left(U_\delta\setminus U_{\delta/2}\right)\approx C_2\cdot\delta
\]
(with some constant $C_2>0$, depending only on $J$) show that there must exist a point 
$y\in U_\delta\setminus U_{\delta/2}$ with the property (a) above,
thus ensuring the sufficiency of the point $\Pi(y)\in J$.

In the closing section of \cite{Sim(2009)} we completed the inductive proof of ergodicity
(with respect to the number of balls $N$) by utilizing Main Lemma 3.5 and
earlier results from the literature.  Actually, a consequence of the
Main Lemma will be that exceptional $J$-manifolds do not exist, and
this will imply the fact that no distinct, open ergodic components can
coexist.

\bigskip

\section{Proof of Ansatz} \label{ansatz}

Finally, in the paper \cite{Sim(2013)} we proved the Boltzmann--Sinai Hypothesis for 
hard ball systems on the $\nu$-torus $\mathbb{R}^\nu/\mathbb{Z}^\nu$ ($\nu\ge 2$) without 
any assumed hypothesis or exceptional model.

As said before, in \cite{Sim(2009)} the Boltzmann-Sinai Hypothesis was proved in full generality
(i.e. without exceptional models), by assuming the Chernov-Sinai Ansatz.

The only missing piece of the whole puzzle is to prove that no
open piece of a singularity manifold can precisely coincide with a
codimension-one manifold desribing the trajectories with a
non-sufficient forward orbit segment corresponding to a fixed symbolic
collision sequence. This is exactly what we prove in our Theorem below.

\subsection{Formulation of Theorem}

Let $U_0\subset \mathbf{M}\setminus\partial\mathbf{M}$ be an open ball, $T>0$, and assume that

\medskip

(a) $S^T(U_0)\cap \partial\mathbf{M}=\emptyset$,

\medskip

(b) $S^T$ is smooth on $U_0$.

\medskip

Next we assume that there is a \emph{codimension-one}, smooth
submanifold $J\subset U_0$ with the property that for every $x\in U_0$
the trajectory segment $S^{[0,T]}x$ is geometrically hyperbolic
(sufficient) if and only if $x\not\in J$. ($J$ is a so called
non-hyperbolicity or degeneracy manifold.) Denote the common symbolic
collision sequence of the orbits $S^{[0,T]}x$ ($x\in U_0$) by
$\Sigma=(e_1,e_2,\dots,e_n)$, listed in the increasing time
order. Let $t_i=t(e_i)$ be the time of the $i$-th collision,
$0<t_1<t_2<\dots<t_n<T$.

Finally we assume that for every phase point $x\in U_0$ the first reflection 
$S^{\tau(x)}x$ \emph{in the past} on the orbit of $x$ is a singular reflection
(i. e. $S^{\tau(x)}x \in \mathcal{SR}^+_0$) if and only if $x$ belongs to a codimension-one,
smooth submanifold $K$ of $U_0$. For the definition of the manifold of singular
reflections $\mathcal{SR}^+_0$ see, for instance, the end of \S1 in \cite{Sim(2009)}.

\begin{thm}
Using all the assumtions and notations above, the submanifolds $J$ and $K$ of $U_0$
do not coincide.
\end{thm}

\end{document}